\documentclass[12pt, reqno]{amsart}
\usepackage[latin1]{inputenc}
\usepackage{amsmath}
\usepackage{amsfonts}
\usepackage{amssymb}
\usepackage[usenames]{color}

\newcommand{\barr}{\bar}

\newcommand{\C}{\mathbb{C}}
\newcommand{\M}{\mathcal{M}}

\newcommand{\CC}{\mathcal{C}}
\newcommand{\HH}{\mathcal{H}}

\newcommand{\X}{\chi}
\newcommand{\T}{\tau}
\newcommand{\f}{\tilde{f}}
\newcommand{\g}{\tilde{g}}

\theoremstyle{definition}
\newtheorem{remark1bb}{Remark}[section]

\theoremstyle{plain}
\newtheorem{theorem1.1a}{Theorem}[section]
\newtheorem{theorem1.2a}[theorem1.1a]{Theorem}
\newtheorem{theorem1.3a}[theorem1.1a]{Theorem}
\newtheorem{theorem1.4a}[theorem1.1a]{Theorem}
\newtheorem{theorem1.5a}[theorem1.1a]{Theorem}
\newtheorem{theorem1.6a}[theorem1.1a]{Theorem}
\theoremstyle{definition}

\theoremstyle{plain}
\newtheorem{lemma2.1a}[theorem1.1a]{Lemma}
\newtheorem{lemma2.2a}[theorem1.1a]{Lemma}
\newtheorem{lemma2.3a}[theorem1.1a]{Lemma}
\newtheorem{lemma2.4a}[theorem1.1a]{Lemma}

\newtheorem{aaLemma2.1}[theorem1.1a]{Lemma}
\newtheorem{lemma2.2aaa}[theorem1.1a]{Lemma}
\newtheorem{aalemma2.3}[theorem1.1a]{Lemma}
\newtheorem{AAAlemma2.4}[theorem1.1a]{Lemma}
\newtheorem{AAAlemma2.5}[theorem1.1a]{Lemma}
\theoremstyle{definition}
\newtheorem{observation1}[theorem1.1a]{Observation}
\newtheorem{exampleA}[theorem1.1a]{Example}

\newtheorem{exampleD}[theorem1.1a]{Example}

\newtheorem{exampleG}[theorem1.1a]{Example}

\newtheorem{exampleJ}[theorem1.1a]{Example}

\newtheorem{exampleO}[theorem1.1a]{Example}

\newtheorem*{theorem1.1}{Theorem 1.1}
\newtheorem*{theorem1.2}{Theorem 1.2}
\newtheorem*{theorem1.3}{Theorem 1.3}
\newtheorem*{theorem1.4}{Theorem 1.4}
\newtheorem*{theorem1.5}{Theorem 1.5}
\newtheorem*{theorem1.6}{Theorem 1.6}
\newtheorem*{lemma2.1}{Lemma 2.1}
\newtheorem*{lemma2.2}{Lemma 2.2}
\newtheorem*{lemma2.3}{Lemma 2.3}
\newtheorem*{lemma2.4}{Lemma 2.4}

\numberwithin{equation}{section}
\begin {document}

\def\1#1{\ov{#1}}
\def\2#1{\widetilde{#1}}
\def\3#1{\mathcal{#1}}
\def\4#1{\widehat{#1}}

\title[Holomorphic Segre Preserving
Maps] {Geometric Properties and Related Results for Holomorphic Segre Preserving
Maps}
\author[R. B. Angle]{R. Blair Angle}
\address{ Department of Mathematics, University of California
at San Diego, La Jolla, CA 92093-0112, USA}
\email{angle@metsci.com }

\maketitle

\begin{abstract}
In this paper, we examine holomorphic Segre preserving maps between the complexifications of real hypersurfaces in $\mathbb{C}^{n+1}$.
In particular, we find
several sufficient conditions  ensuring that  Segre transversality
and total Segre nondegeneracy of the maps must hold.
\end{abstract}

\section{Introduction}

Let $M \subseteq \C^{n+1}$ be a  real analytic hypersurface, with $p \in M$, given locally near $p$ by the real
analytic defining function $\rho(Z,\bar{Z})$.  The \emph{complexification} $\M$ of $M$ is a holomorphic hypersurface of $\C^{2n+2}$ given locally for $(Z,\zeta) \in  \C^{n+1} \times \C^{n+1}$ near $(p,\bar{p})$ by
$\M = \{(Z,\zeta) : \rho(Z,\zeta)=0 \}$.
Let $\Omega \subseteq \C^{n+1}$ be an open neighborhood of $p$ such that $\M$ is defined on $\Omega \, \times \, ^*\Omega$, where
$^*\Omega := \{\barr{Z}:Z \in \Omega\}$.
 Given any $(Z,\zeta) \in \Omega \, \times \, ^*\Omega$,  the \emph{Segre varieties} $\Sigma_Z$ and $\hat{\Sigma}_\zeta$ of $M$ are given by
$\Sigma_Z := \{\zeta\,  \in \, ^*\Omega : \rho(Z,\zeta)= 0\}$ and
$\hat{\Sigma}_\zeta := \{Z\,  \in  \Omega : \rho(Z,\zeta)= 0\}$.
For $(Z',\zeta')$ coordinates on $\C^{n+1} \times \C^{n+1}$,
let $M' \subseteq \mathbb{C}^{n+1}$ be a real analytic hypersurface, with $p' \in M'$,
and denote
its complexification  by $\M'$ and its Segre varieties
 by $\Sigma_{Z'}'$ and $\hat{\Sigma}_{\zeta '}'$. Let $\mathcal{H}:\mathbb{C}^{2n+2}
\rightarrow \mathbb{C}^{2n+2}$  be a holomorphic map defined near
$(p,\bar{p})$ sending $(\M,(p,\bar{p}))$ into $(\M',(p',\bar{p}'))$.  Furthermore, we will assume that for any $(Z,\zeta)
\in \mathcal{M}$, there exists $(Z',\zeta ') \in \mathcal{M}'$ such
that
$\mathcal{H}\big( \{Z\} \times \Sigma_Z \big) \subseteq \{Z'\} \times \Sigma'_{Z'}$ and
 $\mathcal{H}\big( \hat{\Sigma}_\zeta \times\{ \zeta\} \big) \subseteq  \hat{\Sigma}'_{\zeta '} \times \{\zeta'\}$.
  We claim that $\mathcal{H}$, when restricted to $\mathcal{M}$,  is a map of the form
   $$\HH(Z,\zeta)=(H(Z),\widetilde{H}(\zeta)),$$
   where $H,\widetilde{H}: \C^{n+1} \rightarrow \C^{n+1}$.
  This fact was proven for hypersurfaces in \cite{F80}, but it is true for \emph{generic} submanifolds of higher codimension as well (see \cite{Angle}, \cite{A07}).

  We will call such a map a \emph{holomorphic Segre preserving map} (HSPM). Utilizing the notation $\overline{\varphi}(z):= \overline{\varphi(\bar{z})}$, we observe that if  $\widetilde{H}=\overline{H}$, then $H$ is a holomorphic map defined near $p$ sending $(M,p)$ into $(M',p')$. Such maps have been extensively studied. However, HSPMs
are  relatively new and unstudied objects (for related recent work, see  \cite{Angle}, \cite{A07}, and \cite{Z07}).

It can be shown (see, for example, \cite{BER99}) that there exists a holomorphic change of coordinates $Z=(z,w) \in  \C^{n} \times \C$, vanishing at $p$, and an open neighborhood $\Omega$ of 0 such that in these coordinates $M$ is locally given by $\{(z,w) \in \Omega : w=Q(z,\bar{z},\bar{w}) \}$, where $Q(z,\chi,\tau)$ is a holomorphic function defined near 0 in $\C^{n} \times \C^{n} \times \C$
and satisfying
$Q(0,\chi,\tau) \equiv Q(z,0,\tau) \equiv \tau.$
Such coordinates are called \emph{normal coordinates}.

We now express several geometric conditions under the assumption that $M$ is given in normal coordinates (for coordinate-independent definitions see, e.g., \cite{Angle}, \cite{BER99}, \cite{LM07}).
We remind the reader that the \emph{generic rank} of a holomorphic map $\varphi(Z)$ is the largest integer $r$ such that there exists an $r \times r$ minor of the matrix $\left( \frac{\partial \varphi}{\partial Z}\right)$ which is not identically 0. We say that $M$ is of \emph{finite type} at 0 (in the sense of Kohn \cite{Kohn} and Bloom and Graham \cite{BG}) if and only if $Q(z,\X,0) \not\equiv 0$. Otherwise, $M$ is said to be of \emph{infinite type} at 0. In addition, we have:
   
\begin{enumerate}
\item $M$ is \emph{holomorphically nondegenerate} at 0 if and only if  there exists an integer $K$ such that the generic rank of the map
 $(\chi,\tau) \mapsto \big( Q_{z^\alpha}(0,\chi,\tau)\big)_{|\alpha|\leq K}$ is $n+1$.
 \item $M$ is of \emph{class} $\CC$ at 0 if and only if there exists an integer $K$ such that the generic rank of the map
 $\chi \mapsto \big( Q_{z^\alpha}(0,\chi,0)\big)_{|\alpha|\leq K}$ is $n$.
 \item $M$ is \emph{essentially finite} at 0 if and only if there exists an integer $K$ such that  the map
 $\chi \mapsto \big( Q_{z^\alpha}(0,\chi,0)\big)_{|\alpha|\leq K}$ is finite.
 \item $M$ is \emph{finitely nondegenerate} at 0 if and only if  there exists an integer $K$ such that  the matrix
 with rows given by $\big( Q_{\chi_j z^\alpha}(0,0,0)\big)_{|\alpha|\leq K}$, $1 \leq j \leq n$, has rank $n$.
  \end{enumerate}
  In \cite{S96}, Stanton proved that if $M$ is connected and holomorphically
nondegenerate at $0$, then it is holomorphically nondegenerate everywhere. The above definitions thus imply that
finite nondegeneracy at
  0 $\Rightarrow$ essential finiteness at 0 $\Rightarrow$
  class $\CC$ at 0 $\Rightarrow$ holomorphic nondegeneracy. We mention here that the notion of class $\mathcal{C}$ was only
  recently introduced in 2007 by Lamel and Mir in \cite{LM07}.

For the remainder of this paper, we will assume that $M,M' \subseteq \C^{n+1}$ are given in normal coordinates by
$w=Q(z,\X,\T)$ and $w'=Q'(z',\X',\T')$, respectively. We will also assume that
all HSPMs send $(\M,0)$ into $(\M',0)$ and are given in the form
\begin{equation*}
\HH(z,w,\X,\T)= \big(H(z,w),\widetilde{H}(\X,\T)\big) = \big(f(z,w),g(z,w),\f(\X,\T),\g(\X,\T)\big),
\end{equation*}
where $f = (f^1, \ldots, f^n)$ and $\f=(\f^1,\ldots,\f^n)$ are $\C^n$-valued holomorphic functions, and $g$ and $\g$, called the
\emph{transversal} components of $\HH$, are $\C$-valued holomorphic functions. We will write $z=(z_1,\ldots,z_n)$ and $\X=(\X_1,\ldots,\X_n)$
(similarly for $z'$ and $\X'$).

It is easy to see that normality of coordinates implies that $g_w(0)=\g_\T(0)$ (see  \cite{Angle}).
We say that $\HH$ is \emph{Segre transversal} to $\M'$ at 0 if $g_w(0)=\g_\T(0) \neq 0$ (see \cite{Z07}).
$\HH$ is said to be \emph{totally Segre nondegenerate} at 0 if both of the
 following two conditions hold:\
 $\det\big(f_z(z,0)\big) \not\equiv 0$
and $\det\big(\tilde{f}_\chi(\chi,0)\big) \not\equiv 0$.
If only one of these  conditions holds, then $\HH$ is \emph{partially Segre nondegenerate} at 0.
We say that $\HH$ is \emph{transversally null}  if both $g(z,w) \equiv 0$ and $\tilde{g}(\X,\T) \equiv 0$.

\begin{remark1bb}   \label{rem1bb}
It is interesting to note that with HSPMs  it is not necessarily true that
$g \equiv 0 \Leftrightarrow \g \equiv 0$, nor is it necessarily
true that  $\det\big(f_z(z,0)\big) \equiv 0
\Leftrightarrow \det\big(\tilde{f}_\chi(\chi,0)\big) \equiv 0$.  For example, let
$M,M' \subseteq \mathbb{C}^3$ be given by $M=\left\{\text{Im }w=|z_1|^2\right\}$ and $M'=\left\{\text{Im }w'=|z_1'|^2+|z_2'|^2\right\}$.
Thus, $\M=\{w-\T=2iz_1\X_1\}$ and $\M'=\{w'-\T'=2iz_1'\X_1' + 2iz_2'\X_2'\}$.
Define
 $$\mathcal{H}(z,w,\X,\T) =\bigg(z_2+z_1z_2, -z_2, 2wz_2, 2\chi_1, 2\chi_1 +i\T , 0 \bigg).$$
 Observe that $\HH$ sends $(\M,0)$ into $(\M',0)$, but $g \not\equiv 0$, $\tilde{g} \equiv 0$, $\det \big(f_z(z,0)\big) \not\equiv 0$, and $\det \big( \f_\X(\X,0) \big)
 \equiv 0$.
 \end{remark1bb}

If $M'$ is of infinite type at 0, then
given any $\C^n$-valued holomorphic functions $f(z,w)$ and $\f(\X,\T)$ defined near 0 and satisfying $f(0)=\f(0)=0$,
the HSPM  $\mathcal{H}(z,w,\X,\T)=\big(f(z,w),0,\tilde{f}(\X,\T),0\big)$
sends $(\M,0)$ into $(\M',0)$. In fact, if $M$ is of finite type at 0, then these are the
only possible HSPMs (see Lemma \ref{aaLemma2.1aa}). What if we insist that $M'$ be of finite type at 0?
We have the following result which is a generalization of a result of Ebenfelt and Rothschild (\cite{ER05}).

\begin{theorem1.1a} \label{aaThm1}
Let $M,M' \subseteq \mathbb{C}^{n+1}$ be real analytic hypersurfaces, with $M'$  of finite type at $0$, and let $\HH$ be an HSPM
sending $(\M,0)$ into $(\M',0)$. If $\HH$ is totally Segre nondegenerate at 0, then $\HH$ is Segre transversal to $\M'$ at $0$.
\end{theorem1.1a}

Ebenfelt and Rothschild (\cite{ER05}) prove a similar result  for holomorphic maps $H$ sending $M$ into $M'$. The proof of Theorem \ref{aaThm1} follows their proof almost exactly, with only the obvious modifications ($\bar{z}$ and $\bar{w}$ are replaced by $\X$ and $\T$, respectively, and $\bar{f}$ and $\bar{g}$ are replaced by $\f$ and $\g$, respectively).  Thus, we will not present the proof here. (The proof, however, can be found in \cite{Angle}.)

 What if we impose conditions on $M$, rather than $M'$? The next three theorems address this.

\begin{theorem1.2a} \label{aaThm2}
Let $M,M' \subseteq \mathbb{C}^{n+1}$ be real analytic hypersurfaces, with $M$  holomorphically nondegenerate, and let $\mathcal{H}$ be an HSPM sending $(\M,0)$ into $(\M',0)$. Then one of the following two possibilities must hold:
 
\begin{enumerate}
\item $\HH$ is transversally null, and $Q'\left(f(z,w),\tilde{f}(\chi,\tau),0 \right) \equiv 0$. That is, $\mathcal{H}(\mathbb{C}^{2n+2})
\subseteq \mathcal{M}'$.
\item $\HH$ is a biholomorphism near some $p \in \M$.
Furthermore,  $M'$ is holomorphically nondegenerate.
\end{enumerate}
\end{theorem1.2a}

\begin{theorem1.3a} \label{aaThm3}
Let $M,M' \subseteq \mathbb{C}^{n+1}$ be real analytic hypersurfaces, with $M$ of class $\CC$ at $0$,  and let $\mathcal{H}$ be an HSPM
 sending $(\M,0)$ into $(\M',0)$. Then one of the following two possibilities must hold:
  
\begin{enumerate}
\item $\HH$ is transversally null, and $Q'\left(f(z,w),\tilde{f}(\chi,\tau),0 \right) \equiv 0$.
That is, $\mathcal{H}(\mathbb{C}^{2n+2}) \subseteq \mathcal{M}'$.
\item  $\HH$ is Segre transversal to $\M'$ at $0$ and  totally Segre nondegenerate at $0$,
and $\HH$ is a biholomorphism near some $p \in \M$.
 Furthermore, $M'$ is  of class $\mathcal{C}$ at $0$.
\end{enumerate}
\end{theorem1.3a}

\begin{theorem1.4a} \label{aaThm4}
Let $M,M' \subseteq \mathbb{C}^{n+1}$ be real analytic hypersurfaces,
 with  $M$ finitely nondegenerate at $0$,  and let $\mathcal{H}$ be an HSPM sending $(\M,0)$ into $(\M',0)$.
 Then one of the following two possibilities must hold:
  
\begin{enumerate}
\item $\HH$ is transversally null, and $Q'\left(f(z,w),\tilde{f}(\chi,\tau),0\right) \equiv 0$.
That is, $\mathcal{H}(\mathbb{C}^{2n+2}) \subseteq \mathcal{M}'$.
\item $\HH$ is Segre transversal to $\M'$ at $0$, $\det \big(f_z(0)\big) \neq 0$, and
$\det \big( \tilde{f}_\chi(0)\big) \neq 0$. That is, $\mathcal{H}$ is a biholomorphism near $0$.  Furthermore, $M'$ is finitely  nondegenerate at $0$.
\end{enumerate}
\end{theorem1.4a}

Another distinguishing feature of HSPMs is the fact
that since $\tilde{f}$ is not necessarily the complex conjugate of $f$, then the lowest order homogeneous polynomials
 in the Taylor expansions of $\det\big(f_z(z,0)\big)$ and $\det\big(\tilde{f}_\chi(\chi,0)\big)$ may not agree up to
 a constant, assuming neither determinant is identically 0 (see Example  \ref{exJJ}).
 Also, it is possible to have one determinant
  identically 0, and the other not identically 0 (see Remark \ref{rem1bb}).
If $M'$ is finitely nondegenerate at 0, however, we have the following.

\begin{theorem1.5a} \label{aaThm5}
Let $M,M' \subseteq \mathbb{C}^{n+1}$ be real analytic hypersurfaces with $M'$ finitely
nondegenerate at $0$, and let $\mathcal{H}$ be an HSPM,
totally Segre nondegenerate at $0$, sending $(\M,0)$ into $(\M',0)$.
 Assume $\det \big( f_z(z,0)\big) =\sum_{j}p_j(z)$
 and $\det \big( \f_\X(\chi,0)\Big) =\sum_k q_k(\chi)$,
 where $p_l$ and $q_l$ are homogeneous polynomials of degree $l$.
Choose $j_0$ so that $p_{j_0}(z) \not\equiv 0$, but $p_j(z) \equiv 0$ for $j < j_0$. Similarly,
choose $k_0$ so that $q_{k_0}(z) \not\equiv 0$, but $q_k(z) \equiv 0$ for $k < k_0$.
Then $j_0=k_0$, and $p_{j_0}(z)=c\bar{q}_{k_0}(z)$ for some nonzero constant $c$.
\end{theorem1.5a}

This result is essentially sharp. Even if $M'$ is essentially finite at 0, we will see an example in Section \ref{section3}
where the preceding theorem fails. If  $M,M' \subseteq \C^2$, however, then we can be even looser with our hypotheses. We have the following.

\begin{theorem1.6a} \label{aaThm6}
Let $M,M' \subseteq \mathbb{C}^2$ be  real analytic hypersurfaces, with $M'$ of finite type at $0$, and let $\mathcal{H}$ be an HSPM
 sending $(\M,0)$ into $(\M',0)$. Then the following are equivalent:
 
\begin{enumerate}
\item $\HH$ is  totally Segre nondegenerate at $0$.
\item $M$ is of finite type at $0$, and $\HH$ is not transversally null.
\item $M$ is of finite type at $0$, and $\HH$ is Segre transversal to $\M'$ at $0$.
\end{enumerate}
Furthermore, if these conditions hold, then we also have that the order of vanishing
 of $f(z,0)$ equals the order of vanishing of $\tilde{f}(\chi,0)$. That is, if the lowest
 order nonzero term of the Taylor expansion of $f(z,0)$ is $cz^r$ for some nonzero constant $c$ and
 nonnegative integer $r$, then the lowest order nonzero term of the Taylor expansion of $\tilde{f}(\chi,0)$
 is $\tilde{c}\chi ^r$ for some nonzero constant $\tilde{c}$. If, in addition, $M=M'$, then  $\mathcal{H}$ is a biholomorphism near $0$.
\end{theorem1.6a}

The organization of this paper is as follows.  In
Section \ref{section2},  we prove our main results as presented in this section. In Section \ref{section3}
we provide several examples of hypersurfaces and HSPMs between their complexifications.  The purpose of these
examples is to illustrate the necessity of the hypotheses of a particular result or to demonstrate that
a particular result is essentially sharp.

\section{Proofs of main results}   \label{section2}

\subsection{Proof of Theorem \ref{aaThm2} }
Before proving Theorem \ref{aaThm2}, we recall the  definition of holomorphic nondegeneracy.
A hypersurface $M \subseteq \C^{n+1}$ is \emph{holomorphically nondegenerate} at $p$ if there does not exist a germ at $p$ of a nontrivial vector field tangent to $M$ of the form
$L=\sum_{j=1}^{n+1} a_j(Z)\frac{\partial}{\partial Z_j},$
 where $Z=(Z_1,\ldots,Z_{n+1})$ and each $a_j(Z)$ is a holomorphic function.

Now we prove Theorem \ref{aaThm2}.
As $\HH(\M) \subseteq \M'$,
there exists a holomorphic function $A(z,w,\X,\T)$ defined near 0 such that for all $(z,w,\X,\T)$ sufficiently close to 0,
\begin{equation} \label{aaMainEq}
g(z,w)-Q'\big(f(z,w),\tilde{f}(\chi,\tau),\tilde{g}(\chi,\tau)\big)=A(z,w,\chi,\tau)\big(w-Q(z,\chi,\tau)\big) .
\end{equation}
If $A(z,w,\chi,\tau) \equiv 0$, then let $\chi=\tau=0$ in (\ref{aaMainEq}) to see that
$g(z,w) \equiv 0$. Similarly, let $z=w=0$ to see that $\tilde{g}(\chi,\tau) \equiv 0$.
This, in turn, tells us that $Q'\big(f(z,w),\tilde{f}(\chi,\tau),0\big) \equiv 0$. So we now assume that $A(z,w,\chi,\tau) \not\equiv 0$.

If
$\det \big( H_Z(Z) \big)
\equiv 0$, where $Z=(z,w)$, then there exist holomorphic functions $r_j(z,w)$ and $s(z,w)$, not all identically 0, such that
$Lf^1 \equiv \ldots \equiv Lf^n \equiv Lg \equiv 0$, where
\begin{equation}
L := \sum_{j=1}^nr_j(z,w)\frac{\partial}{\partial z_j} + s(z,w)\frac{\partial}{\partial w}.
 \end{equation}
 This is true because there exists a nontrivial vector $V$ of meromorphic functions (the field of fractions of holomorphic functions) annihilated by the matrix $\big(H_Z(Z)\big)$. We write each meromorphic function as the quotient of holomorphic functions, and
 then multiply each component of $V$ by the product of the holomorphic denominators of the meromorphic functions. This gives us $L$.

As $A \not\equiv 0$, we can rewrite (\ref{aaMainEq}) as follows:
\begin{equation} \label{aaMainEq77}
g(z,w)-Q'\big(f(z,w),\tilde{f}(\chi,\tau),\tilde{g}(\chi,\tau)\big)=\widetilde{A}(z,w,\chi,\tau)\big(w-Q(z,\chi,\tau)\big)^m,
\end{equation}
where $\widetilde{A}$ is holomorphic near the origin and not identically zero on $\mathcal{M}$, and $m \in \mathbb{Z}^+$.
Apply $L$ to (\ref{aaMainEq77}) to obtain
\begin{equation*}
0 \equiv L[\widetilde{A}(z,w,\chi,\tau)] \big( w-Q(z,\chi,\tau)\big)^m
\end{equation*}
\begin{equation}
 + \widetilde{A}(z,w,\chi,\tau) m\big(w-Q(z,\chi,\tau) \big)^{m-1}L\left[w-Q(z,\chi,\tau)\right] .
\end{equation}
Dividing both sides by $\big(w-Q(z,\chi,\tau\big))^{m-1}$, we see that $L[w-Q(z,\chi,\tau)] \equiv 0$ on $\M$. This contradicts the holomorphic nondegeneracy of $M$. Thus, our assumption that $\det \big( H_Z(Z) \big) \equiv 0$
 was incorrect. It follows, then, that $H(Z)=H(z,w)$ is a biholomorphism near some point in $\C^{n+1}$.
 A similar argument applies to $\widetilde{H}(\chi,\tau)$. Thus, $\mathcal{H}$ is a biholomorphism near some point in $\C^{2n+2}$.

 In fact, we can find $p_0 \in \M$ such that $\HH$ is a biholomorphism near $p_0$.
 Indeed, for $(\chi,z,w)$ sufficiently close to 0,
$\left\{\big(\chi,\overline{Q}(\chi,z,w)\big)\right\}$ contains an open neighborhood of 0 in $\C^{n+1}$ as $\overline{Q}(\X,0,w) \equiv w$.
 And by assumption, $\det \big(\widetilde{H}_\zeta(\zeta)\big) \not\equiv 0$, where $\zeta=(\X,\T)$.  So choose $\chi_0,z_0,w_0$ sufficiently small so that
\begin{equation} \label{aaEq20}
\det \left(\widetilde{H}_\zeta \big(\chi_0,\overline{Q}(\chi_0,z_0,w_0)\big)\right) \neq 0 .
\end{equation}
As $\det \big(H_Z(Z)\big) \not\equiv 0$, we can choose $(z_0,w_0)$ so that
$\det \big({H}_Z\big(z_0,w_0 \big)\big) \neq 0,$
and also so that (\ref{aaEq20}) holds.

Now we show that $M'$ must be holomorphically nondegenerate. For simplicity, assume $\mathcal{M}$
and $\M'$ are connected and that $\M$ is  defined on $B \times B$, where $B \subseteq \C^{n+1}$ is a ball of sufficiently small radius centered at $0$.
We will also assume, hoping for a contradiction, that $M'$ is holomorphically degenerate.
Thus, it is holomorphically degenerate at all points, in particular at 0.
 Let $$L' := \sum_{j=1}^{n+1} a_j'(Z')\frac{\partial}{\partial Z_j '}$$ be a nontrivial holomorphic vector field defined
 in a neighborhood $U'$ of 0 and tangent to $M'$ in $U'$.
The preceding paragraph implies that there exists $(Z_0',\zeta_0') \in (U' \times \, \hspace{-.035in} ^*U') \cap \M '$ such that $\mathcal{H}^{-1}$ exists near $(Z_0',\zeta_0')$.

      Therefore, near the point $(Z_0',\zeta_0') := \big(H(z_0, w_0), \widetilde{H}(\chi_0,\overline{Q}(\chi_0,z_0,w_0)\big)$, we see that
       if $\rho(Z,\zeta) = w-Q(z,\X,\T)$, then
       $(\rho \circ \HH^{-1})(Z',\zeta ')$ is a defining function for $\M'$.
Thus $L' (\rho \circ \HH^{-1}) = 0$ for $(Z',\zeta') \in \M'$ near $(Z_0',\zeta_0')$.  This then implies that
 $\HH^{-1}_* L' (\rho) = 0$ for $(Z,\zeta) \in \M$ near $\big(z_0, w_0,\chi_0,\overline{Q}(\chi_0,z_0,w_0) \big)$ as
$(Z',\zeta') \in \M'$ if and only if $\HH^{-1}(Z',\zeta') \in \M$.

Note that if $L:=d(z,w)\frac{\partial}{\partial w}+\sum_{j=1}^n c_j(z,w)\frac{\partial}{\partial z_j}$ is tangent to $M$, where $d(z,w)$ and $c_j(z,w)$ are holomorphic functions, then $d(z,w) \equiv 0$. Indeed, apply $L$ to $w-Q(z,\X,\T)$ to see that for all
$(z,w,\X,\T) \in \M$, $d(z,w) - \sum_{j=1}^n c_j(z,w)Q_{z_j}(z,\X,\T) = 0$.  But normality of coordinates implies
that $(z,w,0,w) \in \M$ for all $(z,w)$ sufficiently small.  As $Q_{z_j}(z,0,w) \equiv 0$, it follows that $d(z,w) \equiv 0$.
Thus, we can write $$\HH^{-1}_* L'  = \sum_{j=1}^n a_j(z,w)\frac{\partial}{\partial z_j} .$$ The fact
that this vector field is in terms of $Z$ only, and not $\zeta$, follows from the fact that
 $\mathcal{H}^{-1}(Z',\zeta ')=\big(H^{-1}(Z'),\widetilde{H}^{-1}(\zeta ')\big)$. Write
$$\overline{Q}(\chi,z,w) = \sum_{\alpha}\bar{q}_\alpha(z,w)\chi^\alpha .$$
So we see that
\begin{equation}
\sum_{j=1}^n a_j(z,w)\bar{q}_{\alpha, z_j}(z,w) \equiv 0
\end{equation}
near $(z,w)=(z_0,w_0)$ for all $\alpha$. To complete the proof, we present a lemma which is a slight rewording
of Lemma 11.3.11 in \cite{BER99}.
We remind the reader that each $\bar{q}_\alpha(z,w)$ is defined on the ball $B$.

\begin{lemma2.2aaa} \label{lemma2.2aaa}
Assume there exists $(z_0,w_0) \in B$ and a germ at $(z_0,w_0)$ of a nontrivial $\C^n$-valued
holomorphic function $a(z,w)=\big(a_1(z,w),\ldots,a_n(z,w)\big)$ such that $$\sum_{j=1}^n a_j(z,w)\bar{q}_{\alpha, z_j}(z,w) \equiv 0 ,$$
for all $\alpha$ and all $(z,w)$ near $(z_0,w_0)$. Then there exists a nontrivial
 $\C^n$-valued holomorphic function $b(z,w)=\big(b_1(z,w),\ldots,b_n(z,w)\big)$ defined in $B$ such that
\begin{equation}
\sum_{j=1}^n b_j(z,w)\bar{q}_{\alpha, z_j}(z,w) \equiv 0,
\end{equation}
for all $\alpha$ and all $(z,w) \in B$.
\end{lemma2.2aaa}

To complete the proof of our theorem, we define
\begin{equation}
L := \sum_{j=1}^n b_j(z,w) \frac{\partial}{\partial z_j} ,
\end{equation}
where the $b_j$ are as given in Lemma \ref{lemma2.2aaa}.
Then $L$ is a holomorphic vector field tangent to $M$ on $B$,
 which implies that $M$ is holomorphically degenerate, a contradiction.
$\Box$

\subsection{Proof of Theorem \ref{aaThm3}}
We first present two lemmas.
\begin{aaLemma2.1}  \label{aaLemma2.1aa}
Assume $M$ is of finite type at $0$ and $M'$ is of infinite type at $0$. Then the only HSPMs
 sending $(\M,0)$ into $(\M',0)$ are of the form
\begin{equation}
\HH(z,w,\X,\T) = \big(f(z,w),0,\tilde{f}(\X,\T),0\big),
\end{equation}
for any  $\C^n$-valued holomorphic maps $f(z,w)$ and $\tilde{f}(\X,\T)$ satisfying $f(0)=\tilde{f}(0)=0$.
\end{aaLemma2.1}

\begin{proof}As $\HH$ sends $\M$ into $\M'$,  there exists a holomorphic function $A(z,w,\chi,\tau)$ such that
\begin{equation} \label{mmm}
g(z,w)-Q' \big(f(z,w),\tilde{f}(\chi,\tau),\tilde{g}(\chi,\tau)\big)=A(z,w,\chi,\tau)\big(w-Q(z,\chi,\tau)\big) .
\end{equation}
As normality of coordinates implies that $g(z,0) \equiv \tilde{g}(\X,0) \equiv 0$, and
$M'$ being of infinite type at 0 implies that $Q'(z',\chi',0) \equiv 0$, we substitute $\T=0$ into (\ref{mmm}) to see that
\begin{equation} \label{mmo}
g(z,w)=A(z,w,\chi,0)\big(w-Q(z,\chi,0)\big).
\end{equation}
However, $M$ is of finite type at 0, which implies that $Q(z,\chi,0) \not\equiv 0$. Yet $Q(z,0,0)$ is identically zero.
 Thus, $Q(z,\chi,0)$ must depend on $\chi$. Assume $A(z,w,\chi,0) \not\equiv 0$, and write
 $A(z,w,\chi,0)=\sum_{j=1}^\infty a_j(z,\chi)w^j$. Notice that we
 start our index $j$ at 1 instead of 0.
 Indeed, if we substitute $w=0$ into (\ref{mmo}), we see that $A(z,0,\chi,0) \equiv 0$.
 Let $j_0$ be the smallest positive integer such that $a_{j_0} \not\equiv 0$, but $a_j \equiv 0$ for all $j < j_0$. Then when we
  Taylor expand the right hand side of (\ref{mmo}) in $w$, the lowest order term (i.e., the term with the smallest $w$ exponent) is
   $-a_{j_0}(z,\chi)Q(z,\chi,0)w^{j_0}$. Therefore, the right hand side of (\ref{mmo}) depends on  $\chi$, but the left hand side does not.  This is a contradiction. Thus, our assumption that $A(z,w,\chi,0) \not\equiv 0$ is incorrect, implying that $g(z,w) \equiv 0$.  Similarly $\tilde{g}(\chi,\tau) \equiv 0$.
\end{proof}

\begin{aalemma2.3}  \label{aaLemma2.3aa}
Let $M$, $M'$, and $\HH$ be as in Theorem \ref{aaThm3}. Assume $g(z,w) \not\equiv 0$. Then
$\displaystyle{ \det \big(\tilde{f}_\X(\X,0)\big)} \not\equiv 0$.
Similarly,  $\tilde{g}(\chi,\tau) \not\equiv 0$ implies that
$\displaystyle{ \det \left(f_z(z,0)\right)} \not\equiv 0$.
\end{aalemma2.3}

\begin{proof} We will prove the contrapositive. Assume
 $\displaystyle{ \det \big(\tilde{f}_\X(\X,0)\big)} \equiv 0$. We will show, then,
that $g_{w^k}(z,0) \equiv 0$ for all $k$. This in turn implies that $g(z,w) \equiv 0$. We start with the identity
\begin{equation} \label{aaIdentity}
g\big(z,Q(z,\chi,0)\big)=Q'\Big( f\big(z,Q(z,\chi,0)\big),\tilde{f}(\chi,0),0\Big) .
\end{equation}
Take $\frac{\partial}{\partial \chi}$ of (\ref{aaIdentity}) to get
\begin{displaymath}
\left( g_w\big(z,Q(z,\chi,0)\big)
- \sum_{j=1}^nQ'_{z'_j}\Big(f\big(z,Q(z,\chi,0)\big),\tilde{f}(\chi,0),0\Big)f^j_w\big(z,Q(z,\chi,0)\big) \right)
\end{displaymath}
\begin{equation} \label{aaEq3.47}
 \times \big( Q_{\chi}(z,\chi,0) \big)  = \left( \f_\X
(\chi,0) \right) \Big(Q'_{\chi'} \big(f(z,Q(z,\chi,0)),\tilde{f}(\chi,0),0\big) \Big)  .
\end{equation}
Now assume, looking for a contradiction, that
\begin{equation}  \label{aaEq3.48}
g_w\big(z,Q(z,\chi,0)\big)  - \sum_{j=1}^nQ'_{z_j'}\Big(f\big(z,Q(z,\chi,0)\big),\tilde{f}(\chi,0),0\Big) f_w^j \big(z,Q(z,\chi,0)\big)
\not\equiv 0 .
\end{equation}
 Divide both sides of (\ref{aaEq3.47}) by (\ref{aaEq3.48}) to get
\begin{equation}  \label{aaEq3.49}
\Big( Q_\chi(z,\chi,0)\Big) = \left( \f_\X (\chi,0) \right) \Big(r(z,\chi) \Big)  ,
\end{equation}
where $r(z,\chi)$ is a $\C ^n$-valued holomorphic function, defined near 0 and whenever
 the expression given in (\ref{aaEq3.48}) is not equal to 0. Now  take $\frac{\partial^{|\alpha|}}{\partial z^\alpha}$ of (\ref{aaEq3.49}) to get
\begin{equation}  \label{aaEq1aa}
\Big( Q_{\chi z^\alpha} (z,\chi,0)\Big) = \left( \f_\X (\chi,0) \right) \Big(R_\alpha(z,\chi) \Big)  ,
\end{equation}
where $R_\alpha(z,\chi) := \frac{\partial^{|\alpha|} r}{\partial z^\alpha} (z,\X)$.
  As $M$ is of class $\CC$ at 0, choose  $\alpha_1,\ldots, \alpha_n$ so that the vectors
  $Q_{\chi z^{\alpha_1}} (0,\chi,0), \ldots, Q_{\chi z^{\alpha_n}} (0,\chi,0)$  span $\C^n$ for all $\X \in U \backslash V$, where
$U \subseteq \C^n$ is an open neighborhood of the origin, and $V$ is a proper holomorphic subvariety of $U$.
  Extend (\ref{aaEq1aa}) to the following matrix equation:
  \begin{equation}  \label{bbEq1bb}
\Big( Q_{\chi z^{\alpha_1}} (z,\chi,0), \ldots, Q_{\chi z^{\alpha_n}} (z,\chi,0) \Big) = \left( \f_\X(\chi,0) \right) \Big(R_{\alpha_1}(z,\chi), \ldots, R_{\alpha_n} (z,\chi) \Big)  .
\end{equation}
The matrix on the left hand side of (\ref{bbEq1bb}) is invertible for some $(z,\X)$.
 This contradicts our assumption about $\det \big(\tilde{f}_\chi(\chi,0)\big)$. Thus our assumption  was incorrect. That is, we have
\begin{equation} \label{aaEq30aa}
g_w\big(z,Q(z,\chi,0)\big)  - \sum_{j=1}^nQ'_{z'_j}\Big(f\big(z,Q(z,\chi,0)\big),\tilde{f}(\chi,0),0\Big)f_w^j \big(z,Q(z,\chi,0)\big) \equiv 0.
\end{equation}
Let $\chi=0$ in (\ref{aaEq30aa}) to see that
\begin{equation}
g_w(z,0) \equiv 0 .
\end{equation}

This is the base step of our inductive argument.
Now assume for some $k \geq 1$ we have
\begin{equation}
 g(z,0) \equiv g_w(z,0) \equiv \ldots \equiv g_{w^k}(z,0) \equiv 0,
 \end{equation}
  as well as
\begin{equation}  \label{aaEqttt}
 g_{w^k}\big(z,Q(z,\chi,0)\big) = \sum_{1 \leq |\alpha| \leq k}  Q'_{z'^\alpha}\Big(f\big(z,Q(z,\chi,0)\big),\tilde{f}(\chi,0),0\Big)S^k_\alpha(z,\chi) ,
\end{equation}
where each $S^k_\alpha(z,\chi)$ is the sum of  constant multiples of all possible products of the $f^j_{w^i}\big(z,Q(z,\chi,0)\big)$, with $i \leq k$.
 To complete our induction, we will show that
\begin{equation}
g_{w^{k+1}}\big(z,Q(z,\chi,0)\big) = \sum_{1 \leq |\alpha| \leq k+1}  Q'_{z'^\alpha}\Big(f\big(z,Q(z,\chi,0)\big),\tilde{f}(\chi,0),0\Big)S^{k+1}_\alpha(z,\chi) ,
\end{equation}
where each $S^{k+1}_\alpha(z,\chi)$ is the sum of  constant multiples of all possible products of the
$f^j_{w^i}\big(z,Q(z,\chi,0)\big)$, with $i \leq k+1$.
From this, it immediately follows that
 \begin{equation}
 g_{w^{k+1}}(z,0) \equiv 0  .
 \end{equation}

Take $\frac{\partial}{\partial \chi}$ of (\ref{aaEqttt}). This gives us
\begin{equation*}
\left(g_{w^{k+1}}\big(z,Q(z,\chi,0)\big)  - \sum_{1 \leq |\alpha| \leq k+1}  Q'_{z'^\alpha}\Big(f\big(z,Q(z,\chi,0)\big),\tilde{f}(\chi,0),0\Big)
S^{k+1}_\alpha(z,\chi)\right)
\end{equation*}
\begin{equation}
 \times    \Big( Q_\chi(z,\chi,0) \Big)  = \left( \f_\X (\chi,0) \right) \Big(T(z,\chi) \Big),
\end{equation}
where
$T(z,\X)$ is a $\C ^n$-valued holomorphic function defined near 0.
The induction step is now proved exactly as the base case was proved.
\end{proof}

Now for the proof of Theorem \ref{aaThm3}.
 As $M$ is of class $\mathcal{C}$ at 0, and thus holomorphically nondegenerate,  we see that either the first conclusion of Theorem \ref{aaThm2} holds,
 or we have $g \not\equiv 0$ and $\tilde{g} \not\equiv 0$. Assume the latter. Then Theorem \ref{aaThm2} implies
 that $\HH$ is a biholomorphism near some $p \in \M$.
Lemma \ref{aaLemma2.3aa} implies that  $\det \big(f_z(z,0)\big) \not\equiv 0$ and $\det \big(\tilde{f}_\chi(\chi,0)\big) \not\equiv 0$. From Lemma \ref{aaLemma2.1aa}, we see that $M'$ is of finite type at 0 (since $M$ being of class $\mathcal{C}$ at 0 implies it is of finite type at 0).  Thus, from Theorem \ref{aaThm1} we see that $\HH$ is Segre transversal to $\M'$ at 0.

Now differentiate (\ref{aaIdentity}) with respect to $z_l$, for some $l$:
$$Q_{z_l}(z,\chi,0) = $$
$$\left(g_w\big(z,Q(z,\chi,0)\big)
- \sum_{j=1}^nQ'_{z'_j}\left(f\big(z,Q(z,\chi,0)\big),\tilde{f}(\chi,0),0\right)f^j_w\big(z,Q(z,\chi,0)\big)\right)^{-1} $$
\begin{equation} \label{aba}
\times \left(\sum_{j=1}^n Q'_{z'_j}\left(f\big(z,Q(z,\chi,0)\big),\tilde{f}(\chi,0),0\right)f^j_{z_l}\big(z,Q(z,\chi,0)\big) -
g_{z_l}\big(z,Q(z,\chi,0)\big) \right).
\end{equation}
Notice that the right hand side of (\ref{aba}) is in fact holomorphic near
 0 as $g_w(0) \neq 0$. Thus, we can rewrite (\ref{aba}) in the following way:

$$Q_{z_l}(z,\chi,0) =$$
\begin{equation}  \label{m1}
\frac{ P_l \left( Q'_{z'_j} \left.\left(f\big(z,Q(z,\chi,0)\big),\tilde{f}(\chi,0),0\right) \right|_{1 \leq j \leq n}, H_{z_l}\big(z,Q(z,\chi,0)\big)\right)}
{g_w\big(z,Q(z,\chi,0)\big)
- \sum_{j=1}^nQ'_{z'_j}\left(f\big(z,Q(z,\chi,0)\big),\tilde{f}(\chi,0),0\right)f^j_w\big(z,Q(z,\chi,0)\big)} ,
\end{equation}
where $P_l$ is a  polynomial defined on $\C^n \times
\C^{n+1}$. We can inductively differentiate (\ref{m1}) with respect to $z_k$ for some $k$ to get
$$Q_{z^\alpha}(z,\chi,0) =$$
\begin{equation}    \label{m2}
\frac{ P_\alpha \left( Q'_{z'^\beta} \left. \left( f\big(z,Q(z,\chi,0)\big),\tilde{f}(\chi,0),0 \right)
 \right|_{|\beta|\leq|\alpha|}, H_{z^\beta w^\gamma}\big(z,Q(z,\chi,0)\big)\big|_{|\beta|+|\gamma| \leq |\alpha|}  \right)}
{\Big(g_w\big(z,Q(z,\chi,0)\big)
- \sum_{j=1}^nQ'_{z'_j}\left(f\big(z,Q(z,\chi,0)\big),\tilde{f}(\chi,0),0\right)f^j_w\big(z,Q(z,\chi,0)\big) \Big)^{t_\alpha}} ,
\end{equation}
where $t_\alpha \in \mathbb{Z}^+$ and $P_\alpha$ is a polynomial defined on $\C^{k_\alpha} \times
\C^{l_\alpha}$, for some integers $k_\alpha$ and $l_\alpha$. $\big($Notice that when we differentiate the right
 hand side of (\ref{m2}) with respect to $z_k$, the numerator actually becomes a polynomial function of $Q_{z_k}(z,\chi,0)$ as well, but
we can use (\ref{m1}) to achieve  the form given in (\ref{m2})$\big)$.  Now  set $z=0$ in (\ref{m2}) to get
\begin{equation}  \label{eqR}
Q_{z^\alpha}(0,\chi,0)  = R_\alpha \Big( Q'_{z'^\beta}\big(0,\tilde{f}(\chi,0),0 \big) \big|_{|\beta|\leq|\alpha|}, j_0^{|\alpha|} H \Big) ,
\end{equation}
where $R_\alpha$ is a  rational function
on $\C^{k_\alpha} \times \C^{l_\alpha}$, defined for $\chi$ sufficiently close to 0.
As $M$ is of class $\CC$ at 0, there exist $\alpha_1, \ldots, \alpha_n$  such that
 $Q_{\chi z^{\alpha_1}}(0,\chi,0), \ldots,$ $Q_{\chi z^{\alpha_n}}(0,\chi,0)$  span $\C^n$
for all $\X \in U \backslash V$, where $U \subseteq \C^n$ is an open neighborhood of the origin,
and $V$ is a proper holomorphic subvariety of $U$.
 For convenience, let $k := \text{max}(|\alpha_1|,\ldots,|\alpha_n|)$.
 Using the $R_\alpha$ defined in (\ref{eqR}), we  define a $\C^n$-valued rational
 function $R(y)$ defined (for $y$ sufficiently small) on $\C^{l_k}$, where $l_k$ is the cardinality of the
 set $\{\beta: |\beta| \leq k \}$, such that
\begin{equation*}
\left( \begin{array}{c}
Q_{z^{\alpha_1}}(0,\chi,0), \ldots, Q_{z^{\alpha_n}}(0,\chi,0)
\end{array} \right) = \left(R_{\alpha_j}\left(Q'_{z'^\beta}\big(0,\f(\X,0),0\big)\Big|_{|\beta|\leq|\alpha_j|} , j_0^{|\alpha_j|}H\right)\right)\Bigg|_{1\leq j \leq n}
\end{equation*}
\begin{equation} \label{eqQ}
=R\Big( Q'_{z'^\beta}\big(0,\tilde{f}(\chi,0),0\big) \big|_{|\beta| \leq k} \Big) .
\end{equation}
We differentiate (\ref{eqQ}) with respect to $\chi$, applying the chain rule to get:

$$\left( \begin{array}{c}
Q_{\chi z^{\alpha_1}}(0,\chi,0), \ldots, Q_{\chi z^{\alpha_n}}(0,\chi,0)
\end{array} \right) =$$
\begin{equation}   \label{eqQQ}
\Bigg( \frac{\partial R}{\partial y}\left( Q'_{z'^\beta}\big(0,\tilde{f}(\chi,0),0\big) \big|_{|\beta| \leq k} \right)\Bigg)
\Bigg(\frac{\partial Q'}{\partial z'^\beta \partial \chi'}\big(0,\tilde{f}(\chi,0),0\big)\big|_{|\beta| \leq k} \Bigg)
\Bigg(\frac{\partial \tilde{f}}{\partial \chi}(\chi,0) \Bigg)  .
\end{equation}
As the rank of the matrix on the left hand side of (\ref{eqQQ}) is $n$ for all $\X \in U \backslash V$, we see that the
rank of the $l_k \times n$ matrix $\Big(\frac{\partial Q'}{\partial z'^\beta \partial \chi'}
\big(0,\tilde{f}(\chi,0),0\big)\big|_{|\beta| \leq k} \Big) $ must also be $n$ for all $\X \in U \backslash V$,
 which implies, in particular, that $M'$ is of class $\CC$ at 0.  The proof of Theorem \ref{aaThm3} is complete.
$\Box$

\subsection{Proof of Theorem \ref{aaThm4}}
We now prove Theorem \ref{aaThm4} which follows  immediately from the proof of Theorem \ref{aaThm3}.
 As $M$ is finitely nondegenerate, and thus of class $\CC$, at 0, either the first conclusion
 of Theorem \ref{aaThm3} holds, or the second conclusion holds. Assume the latter. Then
 the fact that $M'$ is finitely nondegenerate at 0 and the fact that $\det\big(\tilde{f}_\chi(0)\big) \neq 0$
follows immediately by letting $\chi=0$ in (\ref{eqQQ}), where we choose $\alpha_1,\ldots,\alpha_n$ so that $Q_{\chi z^{\alpha_1}}(0), \ldots, Q_{\chi z^{\alpha_n}}(0)$
 span $\C^n$.
  Similarly,  $\det\big(f_z(0)\big) \neq 0$.
  $\Box$

\subsection{Proof of Theorem \ref{aaThm5}}
We begin with a lemma.

\begin{AAAlemma2.4}     \label{lleemm}
Let $M$ be  finitely nondegenerate at $0$.
Then for any neighborhood $0 \in U \subseteq \mathbb{C}^n$, there exist $a_1, \ldots, a_n \in U$ such
 that the vectors $\overline{Q}_{z}(a_1,0,0), \ldots, \overline{Q}_{z}(a_n,0,0)$ are linearly independent.
\end{AAAlemma2.4}

\begin{proof} Assume, hoping for a contradiction, that the conclusion is false. Introduce new variables $\eta_1,\ldots, \eta_n \in \mathbb{C}^n$,
and let $\eta =(\eta_1,\ldots, \eta_n)$.  Define the matrix $A(\eta) := \big( \overline{Q}_{z}(\eta_1,0,0),$ $\ldots,
\overline{Q}_{z}(\eta_n,0,0)\big)$. Then there exists an open neighborhood
 $0 \in V \subseteq \mathbb{C}^n$ such that $\det A(\eta) \equiv 0$ for all $\eta_1,\ldots,\eta_n \in V$.
By finite nondegeneracy, assume $\overline{Q}_{z\chi^{\alpha_1}}(0), \ldots,$ $\overline{Q}_{z\chi^{\alpha_n}}(0)$
span $\mathbb{C}^n$. Now take $\displaystyle{\frac{\partial^{|\alpha_1|+\ldots+|\alpha_n|}}{\partial \eta_1^{\alpha_1} \cdots \partial \eta_n^{\alpha_n}}}$ of the identity $\det A(\eta) \equiv 0$ to get
\begin{equation}
\det \bigg(\overline{Q}_{z\chi^{\alpha_1}}(\eta_1,0,0), \ldots, \overline{Q}_{z\chi^{\alpha_n}}(\eta_n,0,0)\bigg) \equiv 0.
\end{equation}
This is clearly false, as finite nondegeneracy implies that this determinant is nonzero when $\eta=0$.
\end{proof}

We now prove Theorem \ref{aaThm5}.
We remind the reader of (\ref{aaEq3.47}) :
\begin{displaymath}
\left(g_w\big(z,Q(z,\chi,0)\big) -  \sum_{j=1}^nQ'_{z'_j}\Big(f\big(z,Q(z,\chi,0)\big),\tilde{f}(\chi,0),0\Big)f^j_w\big(z,Q(z,\chi,0)\big) \right)
\end{displaymath}
\begin{equation}   \label{m3}
  \times \Big( Q_\chi(z,\chi,0) \Big) = \left( \f_\X(\chi,0) \right)
\Big(
Q'_{\chi'} \big( f\left(z,Q(z,\chi,0)\right),\tilde{f}(\chi,0),0\big)  \Big)  .
\end{equation}

Complex conjugation gives
\begin{displaymath}
\left(\bar{g}_\tau\big(\chi,\overline{Q}(\chi,z,0)\big) -
\sum_{j=1}^n \overline{Q}'_{\chi'_j}\Big( \bar{f}\big(\chi,\overline{Q}(\chi,z,0)\big),\bar{\tilde{f}}(z,0),0\Big)\bar{f}^j_\tau\big(\chi,\overline{Q}(\chi,z,0)\big) \right)
\end{displaymath}
\begin{equation}  \label{m4}
\times \Big(
\overline{Q}_z(\chi,z,0) \Big) = \left( \bar{\tilde{f}}_z(z,0) \right)
\Big(
\overline{Q}'_{z'} \big( \bar{f}\big(\chi,\overline{Q}(\chi,z,0)\big),\bar{\tilde{f}}(z,0),0\big)  \Big)  .
\end{equation}
Now we will take $\frac{\partial}{\partial z}$ of the identity
\begin{equation}       \tilde{g}\big(\chi,\overline{Q}(\chi,z,0)\big)=\overline{Q}'\Big(\tilde{f}\big(\chi,\overline{Q}(\chi,z,0)\big),f(z,0),0\Big)    \label{stan2}
\end{equation}
 to get
\begin{displaymath}
\left(\tilde{g}_\tau\big(\chi,\overline{Q}(\chi,z,0)\big) -  \sum_{j=1}^n \overline{Q}'_{\chi'_j}\Big( \tilde{f}\big(\chi,\overline{Q}(\chi,z,0)\big), f(z,0),0\Big)\tilde{f}^j_\tau\big(\chi,\overline{Q}(\chi,z,0)\big) \right)
\end{displaymath}
\begin{equation}
\times \Big( \overline{Q}_z(\chi,z,0) \Big)= \Big( f_z(z,0) \Big)
\left(
\overline{Q}'_{z'} \big( \tilde{f}\big(\chi,\overline{Q}(\chi,z,0)\big), f(z,0),0\big)  \right)  .   \label{m5}
\end{equation}

From Theorem \ref{aaThm1}, as $M'$ is finitely nondegenerate at 0, and thus of finite type at 0, we know that $g_w(0)=\tilde{g}_\tau(0) \neq 0$.
So the Taylor expansions at 0 of both
   \begin{equation}
 \bar{g}_\tau\big(\chi,\overline{Q}(\chi,z,0)\big) - \sum_{j=1}^n \overline{Q}'_{\chi'_j}\Big(\bar{f}\big(\chi,\overline{Q}(\chi,z,0)\big),\bar{\tilde{f}}(z,0),0\Big)\bar{f}^j_\tau\big(\chi,\overline{Q}(\chi,z,0)\big)
 \end{equation} and
  \begin{equation}
 \tilde{g}_\tau\big(\chi,\overline{Q}(\chi,z,0)\big) -  \sum_{j=1}^n \overline{Q}'_{\chi'_j}\Big(\tilde{f}\big(\chi,\overline{Q}(\chi,z,0)\big),f(z,0),0\Big)\tilde{f}^j_\tau\big(\chi,\overline{Q}(\chi,z,0)\big)
    \end{equation}
 have  nonzero constant terms. From Lemma \ref{lleemm}, coupled with our assumptions about $f$ and $\tilde{f}$,
  we can find $a_1,\ldots,a_n$ sufficiently small so that the following hold:
\begin{enumerate}
 \color{white} \item
\color{black}
\begin{enumerate}
    \item[(i)]  $\tilde{g}_\tau(a_j,0) \neq 0$
 \item[(ii)]  $\bar{g}_\tau(a_j,0) \neq 0$
  \item[(iii)]   $\overline{Q}'_z\big(\tilde{f}(a_1,0), 0 ,0\big), \ldots,
 \overline{Q}'_z\big(\tilde{f}(a_n,0),0,0\big)$ span $\mathbb{C}^n$
  \item[(iv)]  $\overline{Q}'_z\big(\bar{f}(a_1,0),0,0\big), \ldots,
   \overline{Q}'_z\big(\bar{f}(a_n,0),0, 0 \big)$ span $\mathbb{C}^n$.
\end{enumerate}
\end{enumerate}

Now define
$$r_j(z):=$$
\begin{equation}
\bar{g}_\tau\big(a_j,\overline{Q}(a_j,z,0)\big) - \sum_{j=1}^n \overline{Q}'_{\chi'_j}
\Big(\bar{f}\big(a_j,\overline{Q}(a_j,z,0)\big),\bar{\tilde{f}}(z,0),0\Big)\bar{f}^j_\tau\big(a_j,\overline{Q}(a_j,z,0)\big),
\end{equation}
 and $$s_j(z):=$$
\begin{equation}
 \tilde{g}_\tau\big(a_j,\overline{Q}(a_j,z,0)\big) -  \sum_{j=1}^n \overline{Q}'_{\chi'_j}\Big(
\tilde{f}\big(a_j,\overline{Q}(a_j,z,0)\big),f(z,0),0\Big)\tilde{f}^j_\tau\big(a_j,\overline{Q}(a_j,z,0)\big).
\end{equation}
Note that  our choice of  $a_j$ implies that $r_j(0) \neq 0$ and $s_j(0) \neq 0$.
 Thus we can form new matrix equations from (\ref{m4}) and (\ref{m5}) by extending the column matrices as follows:
\begin{displaymath}
  \left(
\begin{array}{ccc}
\overline{Q}_{z_1}(a_1,z,0) & \ldots & \overline{Q}_{z_1}(a_n,z,0)  \\
\vdots & \vdots & \vdots  \\
\overline{Q}_{z_n}(a_1,z,0) & \ldots & \overline{Q}_{z_n}(a_n,z,0)
\end{array} \right)
= \left( \begin{array}{ccc}
\bar{\tilde{f}}_{z_1}^1(z,0)  & \ldots & \bar{\tilde{f}}_{z_1}^n(z,0)  \\
\vdots  &   \vdots  &   \vdots  \\
\bar{\tilde{f}}_{z_n}^1(z,0)  & \ldots & \bar{\tilde{f}}_{z_n}^n(z,0)
\end{array} \right) \times
\end{displaymath}
\begin{equation}
\left(
\begin{array}{ccc}
\frac{\overline{Q}'_{z'_1} \Big( \bar{f}\big(a_1,\overline{Q}(a_1,z,0)\big),\bar{\tilde{f}}(z,0),0\Big) }{r_1(z)} & \ldots & \frac{ \overline{Q}'_{z'_1} \Big( \bar{f}\big(a_n,\overline{Q}(a_n,z,0)\big),\bar{\tilde{f}}(z,0),0\Big) }{r_n(z)} \\
\vdots & \vdots  &   \vdots  \\
\frac{ \overline{Q}'_{z'_n} \Big( \bar{f}\big(a_1,\overline{Q}(a_1,z,0)\big),\bar{\tilde{f}}(z,0),0\Big)  }{r_1(z)} & \ldots & \frac{ \overline{Q}'_{z'_n} \Big( \bar{f}\big(a_n,\overline{Q}(a_n,z,0)\big),\bar{\tilde{f}}(z,0),0\Big) }{r_n(z)}
\end{array} \right)  ,  \label{m6}
\end{equation}
and

\begin{displaymath}
  \left(
\begin{array}{ccc}
\overline{Q}_{z_1}(a_1,z,0) & \ldots & \overline{Q}_{z_1}(a_n,z,0)  \\
\vdots & \vdots & \vdots  \\
\overline{Q}_{z_n}(a_1,z,0) & \ldots & \overline{Q}_{z_n}(a_n,z,0)
\end{array} \right)
= \left( \begin{array}{ccc}
f_{z_1}^1(z,0)  & \ldots & f_{z_1}^n(z,0)  \\
\vdots  &   \vdots  &   \vdots  \\
f_{z_n}^1(z,0)  & \ldots & f_{z_n}^n(z,0)
\end{array} \right) \times
\end{displaymath}
\begin{equation}
\left(
\begin{array}{ccc}
\frac{  \overline{Q}'_{z'_1} \Big( \tilde{f}\big(a_1,\overline{Q}(a_1,z,0)\big),f(z,0),0\Big)  }{s_1(z)} & \ldots & \frac{ \overline{Q}'_{z'_1} \Big( \tilde{f}\big(a_n,\overline{Q}(a_n,z,0)\big),f(z,0),0\Big) }{s_n(z)} \\
\vdots & \vdots  &   \vdots  \\
\frac{  \overline{Q}'_{z'_n} \Big( \tilde{f}\big(a_1,\overline{Q}(a_1,z,0)\big),f(z,0),0\Big) }{s_1(z)} & \ldots & \frac{ \overline{Q}'_{z'_n} \Big( \tilde{f}\big(a_n,\overline{Q}(a_n,z,0)\big),f(z,0),0\Big)  }{s_n(z)}
\end{array} \right)   .   \label{m7}
\end{equation}
We set the right hand sides of (\ref{m6}) and (\ref{m7}) equal, and take determinants to see that for $z$ sufficiently small
\begin{equation}
\Big(\det \big( \bar{\tilde{f}}_z(z,0)\big)\Big)\Big( \det C(z) \Big) \left( \prod_{j=1}^ns_j(z) \right) =
\Big(\det \big( f_z(z,0)\big)\Big)\Big( \det D(z) \Big) \left( \prod_{j=1}^nr_j(z) \right)  ,    \label{m8}
\end{equation}

\noindent where $$C(z) := $$
\begin{equation}
\left(
\begin{array}{ccc}
 \overline{Q}'_{z'_1} \Big( \bar{f}\big(a_1,\overline{Q}(a_1,z,0)\big),\bar{\tilde{f}}(z,0),0\Big) & \ldots & \overline{Q}'_{z'_1} \Big( \bar{f}\big(a_n,\overline{Q}(a_n,z,0)\big),\bar{\tilde{f}}(z,0),0\Big) \\
\vdots & \vdots  &   \vdots  \\
\overline{Q}'_{z'_n} \Big( \bar{f}\big(a_1,\overline{Q}(a_1,z,0)\big),\bar{\tilde{f}}(z,0),0\Big) & \ldots &
 \overline{Q}'_{z'_n} \Big( \bar{f}\big(a_n,\overline{Q}(a_n,z,0)\big),\bar{\tilde{f}}(z,0),0\Big)
\end{array} \right)  ,
\end{equation}
and $$D(z):= $$
\begin{equation}
\left(
\begin{array}{ccc}
\overline{Q}'_{z'_1} \Big( \tilde{f}\big(a_1,\overline{Q}(a_1,z,0)\big),f(z,0),0\Big) & \ldots & \overline{Q}'_{z'_1} \Big( \tilde{f}\big(a_n,\overline{Q}(a_n,z,0)\big),f(z,0),0\Big) \\
\vdots & \vdots  &   \vdots  \\
 \overline{Q}'_{z'_n} \Big( \tilde{f}\big(a_1,\overline{Q}(a_1,z,0)\big),f(z,0),0\Big) & \ldots & \overline{Q}'_{z'_n} \Big( \tilde{f}\big(a_n,\overline{Q}(a_n,z,0)\big),f(z,0),0\Big)
\end{array} \right) .
\end{equation}

However, $\det C(0) \neq 0$, $\prod_{j=1}^n s_j(0) \neq 0$,
$\det D(0) \neq 0$, and $ \prod_{j=1}^nr_j(0) \neq 0$.
 Thus, (\ref{m8}) implies that the lowest order homogeneous polynomial in the Taylor expansion of
$\det \big(\bar{\tilde{f}}(z,0)\big)$ is equal (up to a constant) to
the lowest order homogeneous polynomial in the Taylor expansion of
$\det \big(f(z,0)\big)$.  The theorem follows.
$\Box$

\subsection{Proof of Theorem \ref{aaThm6}}

First we show that (i) implies (ii) and (iii).
   Theorem \ref{aaThm1} implies that $\HH$ is Segre transversal to $\M'$ at 0.
The fact that $M$ is of finite type at 0 follows from the identity
\begin{equation} \label{abEq}
g\big(z,Q(z,\chi,0)\big)\equiv Q'\left(f\big(z,Q(z,\chi,0)\big),\tilde{f}(\chi,0),0\right).
 \end{equation}
 Indeed, if $M$ is not of finite type at 0, then $Q(z,\chi,0) \equiv 0$, and equation (\ref{abEq}) becomes
 \begin{equation}
  0 \equiv Q'\big(f(z,0),\tilde{f}(\chi,0),0\big) .
  \end{equation}
  (Normality of coordinates implies that $g(z,0) \equiv 0$.)
  However, this contradicts the fact that $M'$ is of finite type at 0, as $f_z(z,0) \not\equiv 0$ and
  $\f_{\X}(\X,0) \not\equiv 0$.

We note that each of (ii) and (iii)  imply (i).
This follows from the fact that in $\C^2$, being of finite type at 0 is equivalent to being of class $\CC$ at 0.
Thus, (i) follows from Theorem \ref{aaThm3}.

 Now we prove that $f(z,0)$ and $\tilde{f}(\chi,0)$ have the same order of vanishing. As $M$ is of finite type at 0,
we have that $Q(z,\chi,0) \not\equiv 0$.
 Due to normality of coordinates, we can write
\begin{equation}  \label{ass1}
Q(z,\chi,0)=\chi^m\alpha(z) + \chi^{m+1} \beta(z,\chi),
\end{equation}
where $\alpha(z)$ and $\beta(z,\chi)$ are holomorphic functions defined for $z$ and $\chi$ near 0, $m$ is a positive integer, and
\begin{equation}
\alpha(z) = \sum_{j=r}^\infty a_jz^j
\end{equation}
is the Taylor expansion for $\alpha(z)$, with $a_r \neq 0$.
In a similar manner, we can write
\begin{equation}   \label{ass2}
Q(z,\chi,0)=z^n\tilde{\alpha}(\chi) + z^{n+1}\tilde{ \beta}(z,\chi),
\end{equation}
where $\tilde{\alpha}(\chi)$ and $\tilde{\beta}(z,\chi)$ are holomorphic functions, $n$ is a positive integer, and
\begin{equation}  \label{ass2Part2}
\tilde{\alpha}(\chi) = \sum_{j=s}^\infty \tilde{a}_j\chi^j
\end{equation}
is the Taylor expansion for $\tilde{\alpha}(\chi)$, with $\tilde{a}_s \neq 0$.
We make the following observation.

\begin{observation1}  \label{obs1}
$m=n$ and $r=s$.

Indeed, by the reality of our hypersurface,
\begin{equation} \label{ident}
Q\big(z,\chi,\overline{Q}(\chi,z,0)\big) \equiv 0.
 \end{equation}
Due to normality of coordinates, expanding the left hand side of this identity gives us
 \begin{equation}   \label{qqi}
-Q(z,\X,0) \equiv \overline{Q}(\chi,z,0)\big(1 + \chi z \phi(\chi,z)\big),  
\end{equation}
where $\phi(\chi,z)$ is a holomorphic function defined near $\X=z=0$.
Taking the complex conjugate of (\ref{ass1}), it follows that the nonzero term  in the Taylor expansion
 of the right hand side of (\ref{qqi}) with the smallest  $z$ exponent, and then the smallest $\chi$ exponent
given this $z$ exponent, is $\bar{a}_r z^m\chi^r$. It follows then
 that the nonzero term in the Taylor expansion $Q(z,\chi,0)$ with the smallest $z$ exponent, and then the smallest
 $\chi$ exponent given this $z$ exponent, is $-\bar{a}_r z^m\chi^r$.
However, this is equal to $\tilde{\alpha}_s z^n \chi^s$ by (\ref{ass2}). This proves Observation \ref{obs1}.
\end{observation1}

 Now back to the proof of Theorem \ref{aaThm6}.
There exists a   holomorphic function $A(z,\chi,w,$ $\tau)$ defined near 0 such that
\begin{equation}     \label{mm1}
g(z,w)-Q'\big(f(z,w),\tilde{f}(\chi,\tau),\tilde{g}(\chi,\tau)\big) = A(z,\chi,w,\tau)\big(w-Q(z,\chi,\tau)\big).
\end{equation}
Letting $w=\T=0$ in (\ref{mm1}), we have
\begin{equation}     \label{mm2}
Q'\big(f(z,0),\tilde{f}(\chi,0),0\big) = A(z,\chi,0,0)Q(z,\chi,0)  .
\end{equation}

The left hand side of (\ref{mm1}) has a pure $w$ term, as $g_w(0) \neq 0$. Thus,
the right hand side must have one as well, implying that $A(0) \neq 0$. In light of Observation \ref{obs1}, if
the nonzero term  in the Taylor expansion of the right hand side of (\ref{mm2})
with the smallest  $z$ exponent, and then the smallest $\chi$ exponent given this $z$ exponent,
 is $cz^n\chi^s$, for some $c \neq 0$, then the nonzero term with the smallest
 $\X$ exponent, and then the smallest $z$ exponent given this $\X$ exponent, is $d\chi^n z^s$, for some $d \neq 0$.
Thus, the same goes for the left hand side of (\ref{mm2}). As $M'$ is of finite type at 0,
 Observation \ref{obs1} implies that $f(z,0)$ and $\tilde{f}(\chi,0)$ must have the same order of vanishing. In particular, if $M=M'$,
 then $f(z,0)$ and $\tilde{f}(\X,0)$ must have linear terms, implying that $\mathcal{H}$ is a biholomorphism near 0.
 $\Box$

\section{Examples}   \label{section3}
In order to simplify notation in this section, we will drop the convention of using a $'$ to designate
the coordinates of the target hypersurface. The context will render this simplification unambiguous.

In Theorem \ref{aaThm1}, we cannot weaken the hypothesis that $M'$ is of finite type at 0, as the following example illustrates.
\begin{exampleA}  \label{aExA}
Let $M=M' \subseteq \C^2$ be given by
\begin{equation*} \label{ex1.10}
M = \left \{ \text{Im }w=\text{Re }w \left(\frac{\sin |z|^2}{\cos |z|^2 +1}\right)\right\} .
\end{equation*}
Then its complexification is given by
\begin{equation*}    \label{ex1.1}
\M = \{w = \tau e^{iz\chi} \} .
\end{equation*}
Note that $M$ is of infinite type at 0. Let
\begin{equation*}
\HH(z,w,\X,\T) =\left(H(z,w),\widetilde{H}(\X,\T)\right)=\left(\sqrt{2} z,w^2,\sqrt{2} \chi, \tau^2 \right).
 \end{equation*}
This map satisfies the criteria of Theorem \ref{aaThm1}, yet is not
Segre transversal to $\M'$ at 0.
As an aside, notice that $\widetilde{H}=\overline{H}$, and thus the map $H(z,w)=\left(\sqrt{2} z,w^2 \right)$ sends $M$ into itself.
\end{exampleA}

The conclusion of Theorem \ref{aaThm2} does not necessarily hold if $M$ is holomorphically degenerate, as the following
 example illustrates.

\begin{exampleO} \label{Example3.3}
Let $M,M' \subseteq \C^3$ be given by $M= \{\text{Im }w=|z_1|^2\}$ and $M' = \{\text{Im }w=|z_1|^2+|z_2|^2\}$.
Notice that $M$ is holomorphically degenerate and $M'$ is finitely nondegenerate at 0.
Define $\HH(z,w,\X,\T)= (z_1^2,2z_1,w^2 \,,\, 2i\X_1^2,\X_1\T,\T^2)$. Then $\HH$ sends
$(\M,0)$ into $(\M',0)$ and is not transversally null.
Notice, however, that $\HH$ has Jacobian determinant identically zero.
In addition, $\HH$ is neither totally nor partially Segre nondegenerate at 0, nor is $\HH$ Segre transversal to $\M'$ at 0.
\end{exampleO}

The conclusion of Theorem \ref{aaThm2} cannot be improved to that of Theorem \ref{aaThm3}. As an illustration, notice that in Example \ref{aExA},
$\M$ is holomorphically nondegenerate, yet Segre transversality does not hold. Also, Segre nondegeneracy may not hold as the
next example illustrates.

   \begin{exampleD} \label{Example3.6}
Let $M,M' \subseteq \mathbb{C}^3$ be given by $M= \{\text{Im } w =|z_1w|^2 +|z_2|^2 \}$ (defined implicitly)
and $M' = \{\text{Im } w=|z_1|^2 +|z_2|^2 \}$. Notice that $M$ is in fact given in normal coordinates as we have
 \begin{equation*}
\displaystyle{ \M = \left\{ w = \frac{\tau + 2iz_2\chi_2}{1-2iz_1\X_1\T} \right\} }.
\end{equation*}
It can be checked that $M$ is holomorphically nondegenerate and $M'$ is finitely nondegenerate at 0. Define
\begin{equation*}
\mathcal{H}(z,w,\X,\T)= \left(H(z,w),\widetilde{H}(\X,\T)\right) = (z_1w, z_2,w,\chi_1\tau,\chi_2,\tau).
\end{equation*}
 Then $\HH$ sends $(\M,0)$ into $(\M',0)$ and is neither totally nor
partially Segre nondegenerate at 0.
Notice also that $\widetilde{H}=\overline{H}$ and, thus, $H(M) \subseteq M'$.
\end{exampleD}

In Theorem \ref{aaThm5}, we note that  the agreement of the determinants may end after the lowest order homogeneous polynomial.
\begin{exampleG}  \label{exGG}
Let $M=M' \subseteq \mathbb{C}^2$ be the \emph{Lewy hypersurface}
given by $M=\{\text{Im w}=|z|^2\}$. Then $M$ is finitely nondegenerate at 0. Define the following HSPM
which sends $(\M,0)$ into $(\M,0)$:
\begin{equation*}
 \mathcal{H}(z,w,\X,\T) = \Bigg(\frac{2z}{1-2iz},\frac{w}{1-2iz},\frac{1}{2}\chi+\frac{1}{2}\tau,\tau\Bigg) .
\end{equation*}
\end{exampleG}

If $M'$ is finitely degenerate at 0, then the next  example shows that the conclusion of Theorem \ref{aaThm5} does not necessarily hold,
even if $M'$ is essentially finite at 0.

\begin{exampleJ}  \label{exJJ}
Let $M,M' \subseteq \C^6$ be given by
\begin{equation*}
M = \Big\{\text{Im } w=2\text{Re}\big(z_1\bar{z}_2 +z_4^2(\bar{z}_1+\bar{z}_3+\bar{z}_5^2)+
z_2(\bar{z}_3+\bar{z}_5^2)\big) + |z_3|^2+|z_4|^4+|z_5|^4 \Big\},
\end{equation*}
\begin{equation*}
M' = \Big\{\text{Im }w=2\text{Re}(z_3\bar{z}_4^2) + |z_1|^2+|z_2|^2+|z_5|^2 \Big\} .
\end{equation*}
Note that $M$ is finitely nondegenerate at 0, and $M'$ is essentially finite at 0.
The complexifications are given by
$$\mathcal{M}=\left\{ \frac{w-\tau}{2i} =z_1\chi_2+\chi_1z_2 +z_4^2\chi_1+\chi^2_4 z_1 \right.$$
\begin{equation*}
    + z_4^2\chi_3+\chi_4^2z_3+z_4^2\chi_5^2+\chi_4^2z_5^2+z_2\chi_3+\chi_2z_3+
z_2\chi_5^2+\chi_2z_5^2+z_3\chi_3+z_4^2\chi_4^2+z_5^2\chi_5^2 \bigg\},
\end{equation*}
\begin{equation*}
\mathcal{M}'= \bigg\{ \frac{w-\tau}{2i}=z_3\chi_4^2+\chi_3z_4^2+z_1\chi_1+z_2\chi_2+z_5\chi_5 \bigg\}.
\end{equation*}

Let $\mathcal{H}(z,w,\X,\T)=\Big(H(z,w),\widetilde{H}(\X,\T)\Big) =\Big(f(z,w),g(z,w),\f(\X,\T),\g(\X,\T)\Big)$ be defined by
\begin{equation*}
H(z,w)=\big(z_1+z_3,z_4^2+z_2,z_1+z_3+z_4^2,z_5,z_3,w \big),
\end{equation*}
and
\begin{equation*}
\widetilde{H}(\X,\T)= \big( \chi_2,\chi_1+\chi_3+\chi_5^2,\chi_2+\chi_4^2+\chi_5^2,\chi_4,\chi_3,\tau \big).
\end{equation*}
It is easy to show that $\HH(\M) \subseteq \M '$, yet a simple computation proves that
 $\det \big( f_z(z,0)\big) = 2z_4$, and $\det \big( \f_\X(\chi,0)\big) = 2\chi_5$.
\end{exampleJ}

Finally we turn to Theorem \ref{aaThm6}.
  The fact that (i) implies (ii) and (iii) actually holds in $\C^N$ for any $N \geq 2$.
This is clear by examining the first paragraph of the proof of Theorem \ref{aaThm6}.
 However, if we let $r(z)$ denote the lowest order nonzero homogeneous polynomial
 in the Taylor expansion of $\det\big(f_z(z,0)\big)$, then the lowest order nonzero homogeneous polynomial in the Taylor expansion of
$\det\big(\tilde{f}_\chi(\chi,0)\big)$ is not necessarily equal to a constant multiple
of $\bar{r}(\X)$, as Example \ref{exJJ} illustrates.
In $\C^N$ for $N>2$, the other implications, except, of course, (iii) $\Rightarrow$ (ii),
 do not necessarily follow. Example \ref{Example3.3} illustrates that (ii) $\not\Rightarrow$ (i) and
 (ii) $\not\Rightarrow$ (iii). Example \ref{Example3.6} illustrates that (iii) $\not\Rightarrow$ (i).

\end{document}